\documentclass[12pt,thmsa]{article}
\usepackage{amsfonts}
\usepackage{sw20lart}

\input tcilatex2.5
\begin{document}

\title{Some Remarks on Vector-Valued Integration}
\author{V. Kadets \and B. Shumyatskiy \and R. Shvidkoy \and L. Tseytlin and K.
Zheltukhin \\
Department of Mechanics and Mathematics\\
Kharkov State University\\
4 Svobody sq., Kharkov 310077\\
Ukraine\\
\textit{e-mail: }vishnyakova@ilt.kharkov.ua}
\date{February, 1999}
\maketitle

\begin{abstract}
The article presents a new method of integration of functions with values in
Banach spaces. This integral and related notions prove to be a useful tool
in the study of Banach space geomtry.
\end{abstract}

\section{Introduction}

There are many definitions of integral for Banach-space-valued functions and
their corresponding classes of integrable functions. The easiest one is the
Riemann integral. This definition looks exactly like the one for real-valued
functions.

\begin{definition}
Let $f:[0;1]\rightarrow X$ be a bounded function. This function is said to
be \emph{Riemann integrable} if there exists an $x\in X$ (called $f$'s \emph{%
Riemann integral}) such that for any $\epsilon >0$ there is a $\delta >0$
such that for any partition of the segment $[0;1]$ into a finite number of
intervals $\{\Delta _i\}_{i=1}^N$ with $\max_i\left| \Delta _i\right|
<\delta $ and any choice of sampling points $t_i\in \Delta _i$ the
corresponding \emph{Riemann integral sum} $S_R(f,\{\Delta
_i\},\{t_i\})=\sum_{i=1}^Nf(t_i)\left| \Delta _i\right| $ is close to $x$,
that is, $\left\| S_R(f,\{\Delta _i\},\{t_i\})-x\right\| <\epsilon $.
\end{definition}

The most commonly used definition of a vector-valued integral is, however,
that of Bochner integral.

\begin{definition}
Let $(\Omega ,\Sigma ,\mu )$ be a measure space and $f:\Omega \rightarrow X$
be a function. The function $f$ is called \emph{Bochner-integrable} if there
exists a sequence of simple (measurable finite-valued) functions $f_n$ such
that $\int_\Omega \left\| f-f_n\right\| d\mu \rightarrow 0,$ as $%
n\rightarrow \infty $. The \emph{Bochner integral} of $f$ is then defined as 
$\int_\Omega fd\mu =\lim_{n\rightarrow \infty }\int_\Omega f_nd\mu $, where
the integral for simple functions is defined in the obvious way: if $%
f_n=\sum x_i\chi _{A_i}$, then $\int_\Omega f_nd\mu =\sum x_i\mu (A_i)$. It
is easy to see that such a limit indeed exists and does not depend on the
choice of a sequence $f_n$ approximating the given function $f$.
\end{definition}

Although the Bochner integrability is a direct generalization of the
Lebesgue one and has many properties of the latter, there are Riemann
integrable functions that are not Bochner integrable.

\emph{Example. }\label{ex-Riemann-non-Bochner}Consider the space $l_2([0;1])$%
. It consists of all functions $f:[0;1]\rightarrow \Bbb{R}$ such that $%
\sum_{t\in [0;1]}|f(t)|^2<\infty $ (it follows that these functions take
non-zero values on countable subsets of $[0;1]$). The norm on $l_2([0;1])$
is given by $||f||=(\sum_{t\in [0;1]}|f(t)|^2)^{1/2}$. It is easy to see
that $l_2([0;1])$ is a non-separable Hilbert space. Its orts are given by $%
e_t=\chi _{\{t\}}$.

Now consider a function $f:[0;1]\rightarrow l_2([0;1])$: $f(t)=e_t$. This
function is not measurable, since it is not even ``almost
separable-valued''. Thus it is not Bochner-integrable. However, it is
Riemann-integrable. To see this, note that for any partition $\Pi =\{\Delta
_i\}$ with any sample points $T=\{t_i\}$, 
\begin{eqnarray*}
||S_R(f,\Pi ,T)|| &=&||\sum f(t_i)\left| \Delta _i\right| ||=||\sum
e_{t_i}\left| \Delta _i\right| ||= \\
\ &=&(\sum \left| \Delta _i\right| ^2)^{1/2}\leq (\sum d(\Pi )\cdot \left|
\Delta _i\right| )^{1/2}= \\
\ &=&\sqrt{d(\Pi )}(\sum_i\left| \Delta _i\right| ))^{1/2}=\sqrt{d(\Pi )},
\end{eqnarray*}
where $d(\Pi )$ denotes $\max \{\left| \Delta _i\right| \}$. Now take an
arbitrary $\epsilon >0$. Fix a partition $\Pi $ with $d(\Pi )<\epsilon $.
Then for any $\Pi ^{\prime }\succ \Pi $ and any $T^{\prime }$ we have $%
||S_R(f,\Pi ^{\prime },T^{\prime })||\leq \sqrt{d(\Pi ^{\prime })}\leq \sqrt{%
d(\Pi )}\leq \epsilon $. Thus $f$ is Riemann-integrable and the integral
equals $0$.

\bigskip\ 

This gap may be covered by the idea of weak integration.

\begin{definition}
Let $(\Omega ,\Sigma ,\mu )$ be a measure space and $f:\Omega \rightarrow X$
be a function. This function is called \emph{weakly integrable} if for any
functional $F\in X^{*}$ the scalar function $F\circ f$ is
Lebesgue-integrable.

The function $f$ is called \emph{Pettis-integrable} if for any $A\in \Sigma $
there exists an $x\in X$ such that $F(x)=\int_AF\circ fd\mu $, whenever $%
F\in X^{*}$. The point $x$ is then called $f$'s \emph{Pettis integral} over $%
A$.
\end{definition}

Pettis integrability is useful enough and covers both Riemann and Bochner
types of integration, but it, in turn, has its own drawbacks. The definition
of the Pettis integral does not use any ``simple'' approximations, that is
why some of the usual Lebesgue-integration theorems do not hold for the
Pettis integral. For example, there exists a Pettis-integrable function $%
f:[0;1]\rightarrow X$, for which the ``antiderivative'' $F(t)=\int_0^tf(\tau
)d\mu $ is not differentiable (\cite{Kadets}); the space of all
Pettis-integrable functions is not complete, etc. We refer the reader to
classical texts \cite{Tal} and \cite{Diest-Uhl} for the detailed treatment
of this theory.

In the present paper we study a new definition of integrability, introduced
by two of the authors in \cite{Kadets-Tseytlin} and named the RL
(Riemann-Lebesgue) integrability, which covers both Riemann and Bochner
integrals (theorems \ref{Riemann} and \ref{th-Bochner-RL}), but is less
general than the Pettis one. As will be shown below, this notion is a very
natural and convenient one. The properties of the RL integration depend on
the properties of the space $X$, which makes this concept a valuable tool
for the study of Banach space structure (especially in the non-separable
case).

\section{RL-Integral Sums and RL-Integral}

Let $(\Omega ,\Sigma ,\mu )$ be a space with finite measure $\mu $ and $X$
be a Banach space.

\begin{definition}
Let $f:\Omega \rightarrow \Bbb{R}^{+}$ be a function. We call the value $%
\overline{\int }fd\mu =\inf \{\int gd\mu :\,g(t)\geq f(t)\,\forall t;\,g\,%
\text{is Lebesgue integrable}\}$ $f$'s \emph{upper Lebesgue integral}. In
particular, $\overline{\int }fd\mu =+\infty $, if $f$ has no a
Lebesgue-integrable majorant.
\end{definition}

Using the notion of the upper Lebesgue integral we introduce the \emph{upper}%
-$L_1$ ($\overline{L_1}$) \emph{space. }

\begin{definition}
The space $\overline{L_1}(\Omega ,\Sigma ,\mu ,X)$ is defined as the space
of all functions $f:\Omega \rightarrow X$ such that $\overline{\int }\left\|
f\right\| d\mu <+\infty $. The norm on this space is given by $\left\|
f\right\| =\overline{\int }\left\| f(t)\right\| d\mu (t)$, and is called the 
\emph{RL-norm} ($\left\| \cdot \right\| _{RL}$).
\end{definition}

It is easy to see that any bounded function as well as any
Bochner-integrable function belongs to $\overline{L_1}$. Using usual
classical methods, one can also prove that $\overline{L_1}$ is a Banach
space.

Let us now introduce the notion of RL integral sums, first defined in \cite
{Kadets-Tseytlin}.

\begin{definition}
Let $f:\Omega \rightarrow X$ be a function. Let $\Pi =\{\Delta
_i\}_{i=1}^\infty $ be a partition of $\Omega $ into a countable number of
measurable subsets. Let $T=\{t_i\}_{i=1}^\infty $ be the set of sampling
points for $\Pi $, i.e. $t_i\in \Delta _i$. We can construct a formal series 
$S(f,\Pi ,T)=\sum_{i=1}^\infty f(t_i)\mu (\Delta _i)$. This series is called
the \emph{(absolute) RL integral sum} of $f$ with respect to $\Pi $ and $T$,
provided it is absolutely convergent.
\end{definition}

It seems natural to order the set of partitions in the following way. We
will say that the partition $\Pi _1=\{D_i\}$ \textit{follows} partition $\Pi
_2=\{E_k\}$, or $\Pi _1$ \textit{is inscribed into} $\Pi _2$ ($\Pi _1\succ
\Pi _2$), whenever $\Pi _1$ is a finer partition, that is, the set of
indices $\Bbb{N}$ can be broken into disjoint subsets $I_k$, $k\in \Bbb{N}$
, $\bigcup_{k=1}^\infty I_k=\Bbb{N}$ such that $E_k=\bigcup_{i\in I_k}D_i$.

\begin{lemma}
Let $f\in \overline{L_1}(\Omega ,\Sigma ,\mu ,X)$. Then there exists a
partition $\Pi $ of $\Omega $ so that for any finer partition $\Gamma \succ
\Pi $ and any set of sampling points $T$ the series $S(f,\Gamma ,T)$ is
absolutely convergent.
\end{lemma}

\emph{Proof. }Since $f\in \overline{L_1}(\Omega ,\Sigma ,\mu ,X)$, there
exists a Lebesgue-integrable function $g:\Omega \rightarrow \Bbb{R}^{+}$
that dominates the norm $f(t)$ pointwise. Since $g$ is Lebesgue-integrable,
we can choose an $\epsilon >0$ so small that the upper Lebesgue integral sum 
$\sum_{i=1}^\infty i\epsilon \cdot \mu (g^{-1}([(i-1)\epsilon ;\,i\epsilon
]))$ is absolutely convergent. Denoting $\Delta _i=g^{-1}([(i-1)\epsilon
;\,i\epsilon ])$ we obtain a partition $\Pi $ for which the series $%
\sum_i\sup_{t\in \Delta _i}g(t)\cdot \mu (\Delta _i)$ is absolutely
convergent.

Now if $\Gamma \succ \Pi $, then we can write $\Gamma =\{\Delta _{ij}\}$,
where $\bigcup_j\Delta _{ij}=\Delta _i$. Let $T=\{t_{ij}\}$ be a set of
sampling points for $\Gamma $. Then the integral sum $S(f,\Gamma
,T)=\sum_{ij}f(t_{ij})\mu (\Delta _{ij})$ is dominated in norm by the series 
$\sum_{ij}\sup_{t\in \Delta _i}g(t)\mu (\Delta _{ij})$, which is convergent.
Thus $S(f,\Gamma ,T)$ is absolutely convergent. $\Box $

\bigskip\ Now we can introduce the definition of the RL integral.

\begin{definition}
A function $f:\Omega \rightarrow X$ is called \textit{Riemann-Lebesgue} (%
\textit{RL}) \textit{integrable} over a measurable set $A\subset \Omega $ if
there exists a point $x\in X$ such that for any $\epsilon >0$ there is a
partition $\Pi $ of $A$ such that for any finer partition $\Gamma \succ \Pi $
with any set of sampling points $T$ we have $||S(f,\Gamma ,T)-x||<\epsilon $
and the sum $S(f,\Gamma ,T)$ converges absolutely. This point $x$ is called
then the Riemann--Lebesgue integral of $f$ and denoted, as usual, by $%
\int_Af(t)dt$.
\end{definition}

Let us point out several simple properties of the RL integral:

\begin{enumerate}
\item  RL integral is additive: if functions $f$ and $g$ are RL-integrable,
then $f+g$ is RL-integrable too and $\int (f+g)d\mu =\int fd\mu +\int gd\mu $%
.

\item  If a function $f$ is RL-integrable over a set $A$, it is also
RL-integrable over any measurable subset of $A$.

\item  If $f:\Omega \rightarrow X$ is an RL-integrable function, then $f\in 
\overline{L_1}(\Omega ,\Sigma ,\mu ,X)$.

\item  As a function of set, RL integral is a countably-additive measure of
bounded variation.

\item  If $f:\Omega \rightarrow X$ is an RL-integrable function and $%
T:X\rightarrow Y$ is a continuous linear operator, then the composition $%
Tf:\Omega \rightarrow Y$ is an RL-integrable function and we have $T\int
fd\mu =\int Tfd\mu $.
\end{enumerate}

\begin{definition}
The $RL_1(\Omega ,\Sigma ,\mu ,X)$ \emph{space} is the linear subspace of $%
\overline{L_1}(\Omega ,\Sigma ,\mu ,X)$ consisting of all those functions
that are RL-integrable.
\end{definition}

$RL_1$ is a closed subspace of $\overline{L_1}$ and as such is a Banach
space.

The following theorem is a useful sufficient condition for a function to be
RL-integrable. We will use it later to prove the RL-integrability of the
Bochner and Riemann integrable functions.

\begin{theorem}
\label{th-seq-RL-integrable}Let $f\in \overline{L_1}(\Omega ,\Sigma ,\mu ,X)$
be a function and $g\in L_1(\Omega ,\Sigma ,\mu )$ be its integrable
majorant. Let also $\{\Pi _n\}$ be a sequence of partitions, for each of
which the upper integral sum of $g$ is absolutely convergent. Assume that
for any choice of sampling points $T_n$ for $\Pi _n$, the sequence $S(f,\Pi
_n,T_n)$ converges to a certain element $x\in X$. Then $f$ is RL-integrable
and $x=\int fd\mu $.
\end{theorem}

Let us prove some auxiliary lemmas first.

\begin{lemma}
\label{lem-conv-subset}Let $K_1,\ldots ,K_N$ be finite sets in a Banach
space $X$. Then 
\begin{equation}
\limfunc{conv}(K_1)+\ldots +\limfunc{conv}(K_N)\subset \limfunc{conv}%
(K_1+\ldots +K_N)  \label{eqn-conv-subset}
\end{equation}
\end{lemma}

\emph{Proof. }Since the sets $K_i$ are finite and their number is finite,
both sides of (\ref{eqn-conv-subset}) are convex compacts. Due to the
Krein--Milman theorem, it is sufficient to show that the extreme points of
the left-hand side of (\ref{eqn-conv-subset}) belong to the right-hand one.

Denote the left-hand side of (\ref{eqn-conv-subset}) by $A$. Let $x$ be an
extreme point of $A$. Obviously, $x=x_1+\ldots +x_N$, where $x_i\in \limfunc{%
conv}(K_i)$, $1\leq i\leq N$. Note that each $x_i$ must be an extreme point
of $\limfunc{conv}(K_i)$. Indeed, assume that one of the $x_i$'s is not
extreme. Without loss of generality we can assume that $i=1$. So $%
x_1=(y_1+z_1)/2$, where $y_1,z_1\in \limfunc{conv}(K_1)$. But then $%
x=(y+z)/2 $, where $y=y_1+x_2+\ldots +x_n$ and $z=z_1+x_2+\ldots +x_n$, $%
y,z\in A$, which is impossible, since $x$ is an extreme point of $A$.

Now note that an extreme point of $\limfunc{conv}(K_i)$ must belong to $K_i$
itself. Thus $x_i\in K_i$, and therefore $x=x_1+\ldots +x_N\in K_1+\ldots
+K_N\subset \limfunc{conv}(K_1+\ldots +K_N)$, which proves the lemma. $\Box $

\begin{lemma}
\label{lem-finite-finer-subset}Let $f:\Omega \rightarrow X$ be a function.
Let $\Pi $ be a finite partition of $\Omega $. Then for any finer finite
partition $\Gamma \succ \Pi $ the following inclusion holds: 
\begin{eqnarray*}
\{S(f,\Gamma ,T) &:&\,T\text{ any choice of sampling points}\}\subset \\
\ &\subset &\limfunc{conv}\{S(f,\Pi ,T):\,T\text{ any choice of sampling
points}\}.
\end{eqnarray*}
\end{lemma}

\emph{Proof. }Let $\Pi =\{\Delta _i\}_{i=1}^N$ and let $\Gamma =\{\Delta
_{ij}\}_{i=1}^N\,_{j=1}^{M_i}$ where $\bigcup_{j=1}^{M_i}\Delta _{ij}=\Delta
_i$. Consider arbitrary sampling points $T=\{t_{ij}\}$, $t_{ij}\in \Delta
_{ij}$ and the integral sum $S(f,\Gamma ,T)=\sum_{ij}f(t_{ij})\mu (\Delta
_{ij})$. Denote by $S_i$ the part of this sum where $\Delta _{ij}\subset
\Delta _i$, i.e. $S_i=\sum_jf(t_{ij})\mu (\Delta _{ij})$. Then $S(f,\Gamma
,T)=S_1+\ldots +S_N$. Denote also $K_i=\{f(t_{ij})\mu (\Delta
_i):\,j=1,\ldots .M_i\}$ for each $i=1,\ldots ,N$.

Note that 
\[
S_i=\sum_{j=1}^{M_i}f(t_{ij})\mu (\Delta _i)\cdot \frac{\mu (\Delta _{ij})}{%
\mu (\Delta _i)}, 
\]
where $\sum_j\frac{\mu (\Delta _{ij})}{\mu (\Delta _i)}=1$. Thus $S_i$ is a
convex combination of elements of the form $f(t_{ij})\mu (\Delta _i)$, i.e. $%
S_i\in \limfunc{conv}(K_i)$. Due to lemma~\ref{lem-conv-subset}, $S(f,\Gamma
,T)\in \limfunc{conv}(K_1+\ldots +K_N)$. But each element of $K_1+\ldots
+K_N $ has the form $f(t_{1j_1})\mu (\Delta _1)+\ldots +f(t_{Nj_N})\mu
(\Delta _N) $, where $t_{ij_i}\in \Delta _{ij_i}\subset \Delta _i$. Thus
each element of $K_1+\ldots +K_N$ is an integral sum of the form $S(f,\Pi
,T^{\prime })$ for certain choice of sampling points $T^{\prime }$. This
proves the lemma. $\Box $

\begin{lemma}
\label{lem-infinite-finer-subset}Let $f\in \overline{L}_1(\Omega ,\Sigma
,\mu ,X)$ be a function and $g\in L_1(\Omega ,\Sigma ,\mu )$ be its
integrable majorant. Let $\Pi $ be a partition of $\Omega $ for which the
upper integral sum of $g$ is absolutely convergent. Then for any finer
partition $\Gamma \succ \Pi $ the following inclusion holds: 
\begin{eqnarray*}
\{S(f,\Gamma ,T):\,T\text{ any choice of sampling points}\}\subset \\
\ \subset \overline{\limfunc{conv}}\{S(f,\Pi ,T):\,T\text{ any choice of
sampling points}\}.
\end{eqnarray*}
\end{lemma}

\emph{Proof.} The lemma can be reduced to the previous one by approximating
infinite series with their partial sums. $\Box $

\bigskip Now we are in a position to prove theorem~\ref{th-seq-RL-integrable}%
.

\emph{Proof of theorem \ref{th-seq-RL-integrable}}$.$

Note that under the conditions of the theorem, the following fact is true:
for any $\epsilon >0$ there exists a partition $\Pi _n$ such that for any
choice of sampling points $T_n$ for it, $\left\| S(f,\Pi _n,T_n)-x\right\|
<\epsilon $. Indeed, if this were not true, then there would exist an $%
\epsilon >0$ for which we could find a subsequence $\{\Pi
_{n_k}\}_{k=1}^\infty $ with sampling points $T_{n_k}$ so that $\left\|
S(f,\Pi _{n_k},T_{n_k})-x\right\| \geq \epsilon $. Then choosing arbitrary
sets of sampling points $T_n$ for partitions $\Pi _n$, where $n\notin
\{n_k\} $, we would have a sequence $(\Pi _n,T_n)$ for which $S(f,\Pi
_n,T_n) $ does not converge to $x$.

Now fix an $\epsilon >0$. Take a partition $\Pi _n$ for which $\left\|
S(f,\Pi _n,T_n)-x\right\| <\epsilon $ for any choice of sampling points $T_n$%
. Then, by lemma~\ref{lem-infinite-finer-subset}, for any $\Gamma \succ \Pi
_n$ and any choice of sampling points $T$ for $\Gamma $ we have $\left\|
S(f,\Gamma ,T)-x\right\| <\epsilon $, which proves that $f$ is RL-integrable
and the RL-integral of $f$ equals $x$. $\Box $

\bigskip This theorem in hand we can show that the notion of
RL-integrability of functions is more general than some known ones.

\begin{theorem}
\label{Riemann}Let a bounded function $f:[0;1]\rightarrow X$ be Riemann
integrable. Then it is RL-integrable and its RL integral equals its Riemann
integral.
\end{theorem}

\emph{Proof. }Indeed, since $f$ is bounded, it belongs to $\overline{L}_1$.
The constant function that bounds $f$ can serve as the dominant $g$ from the
theorem~\ref{th-seq-RL-integrable}. Its upper Lebesgue integral sum for any
partition is obviously convergent. Let $\Pi _n$ be a partition of the
segment $[0;1]$ into $2^n$ equal subsegments. By the definition of the
Riemann integral, $\{\Pi _n\}$ is a sequence of partitions that satisfies
the conditions of theorem~\ref{th-seq-RL-integrable}. Therefore it is
RL-integrable. $\Box $

\begin{theorem}
\label{th-Bochner-RL}Let a function $f:\Omega \rightarrow X$ be
Bochner-integrable. Then it is RL-integrable and its RL-integral equals its
Bochner integral.
\end{theorem}

\emph{Proof.} We use theorem~\ref{th-seq-RL-integrable}. Since $f$ is
Bochner integrable, the function $g(t)=\Vert f(t)\Vert $ is its integrable
majorant and hence $f$ belongs to $\overline{L_1}(\Omega ,\Sigma ,\mu ,X)$.
Put $x=(Bochner)\int fd\mu $. Let us fix a sequence of positive numbers $%
\epsilon _n\rightarrow 0$ and construct the partition $\Pi _n$ required in
theorem \ref{th-seq-RL-integrable}.

Due to $f$'s Bochner integrability, it is measurable and therefore almost
separable-valued. Thus there exists a subset $\Omega ^{\prime }\subset
\Omega $, $\mu (\Omega ^{\prime })=\mu (\Omega )$, such that $f(\Omega
^{\prime })$ is separable. Let us cover the separable set $f(\Omega ^{\prime
})$ with countable number of disjoint sets $A_i$ so that $\limfunc{diam}%
(A_i)<\epsilon _n$. To do this we first fix a countable dense set $%
\{x_i\}\subset f(\Omega ^{\prime })$, then consider balls $B_i=\{x\in
X:\,\Vert x-x_i\Vert <\epsilon _n\}$, and then define $A_i$'s as follows: $%
A_1=B_1$, $A_2=B_2\setminus A_1$, $A_3=B_3\setminus (A_1\cup A_2)$, and so
on. It is clear from the construction that the sets $f^{-1}(A_i)$ are
measurable. Therefore $\{f^{-1}(A_i)\}$ is a partition of $\Omega ^{\prime }$%
. By adding a negligible set $\Omega \setminus \Omega ^{\prime }$ (which
will have no effect on the integral sums) to the partition we obtain a
partition of the entire $\Omega $. Denote this partition by $\Pi _n=\{\Delta
_i\}$. Obviously, it has the following property: 
\begin{equation}
\text{if }\mu (\Delta _i)\neq 0\text{, then }\limfunc{diam}(f(\Delta
_i))<\epsilon _n.  \label{eqn-Bochner-RL-1}
\end{equation}

Let us prove that $\Pi _n$ is required. To this end, suppose $T=\{t_i\}$ is
a set of sampling points for $\Pi _n$. Then, taking into account (\ref
{eqn-Bochner-RL-1}), we estimate: 
\begin{eqnarray*}
\left\| f(t_i)\mu (\Delta _i)-\int_{\Delta _i}fd{\mu }\right\| &=&\left\| {%
\int_{\Delta _i}}f(t_i)d\mu {-}\int_{\Delta _i}fd{\mu }\right\| \leq \\
&\leq &\left\| \int_{\Delta _i}f(t_i)-f(t)d\mu \right\| \leq \epsilon _n\mu
(\Delta _i).
\end{eqnarray*}
Therefore 
\begin{eqnarray*}
\left\| S(f,\Pi _n,T)-x\right\| &=&\left\| \sum_{i=1}^\infty f(t_i)\mu
(\Delta _i)-\sum_{i=1}^\infty \int_{\Delta _i}fd\mu \right\| \leq \\
&\leq &\left\| \sum_{i=1}^\infty \left( f(t_i)\mu (\Delta _i)-\int_{\Delta
_i}fd\mu \right) \right\| \leq \epsilon _n\sum_{i=1}^\infty \mu (\Delta
_i)=\epsilon _n.
\end{eqnarray*}

Finally, to show that for the partition $\Pi _n$ the upper Lebesgue integral
sum of the majorant $g(t)=\Vert f(t)\Vert $ is convergent, note that the
series $\sum \sup_{\Delta _i}g\cdot \mu (\Delta _i)$ is dominated by the
absolutely convergent series $\sum \int_{\Delta _i}(g+\epsilon _n)d\mu $.
This finishes the proof. $\Box $

\bigskip Now let us mention two properties of the RL integral proved in~\cite
{Kadets-Tseytlin}.

\begin{enumerate}
\item  An RL-integrable function is Pettis-integrable and the values of both
integrals coincide.

\item  If the space $X$ is separable, then RL-integrability is equivalent to
Bochner integrability.
\end{enumerate}

\section{$RL_1(X,\mu )$ Space}

Let us investigate the properties of the space $RL_1(X,\mu )$. As a closed
subspace of the space $\overline{L_1}(X,\mu )$, it is a Banach space. Note
that the space $L_1(X,\mu )$ of Bochner-integrable $X$-valued functions is a
closed subspace of $RL_1(X,\mu )$. It seems natural to ask: under what
condition on the space $X$, $L_1(X,\mu )$ is complemented in $RL_1(X,\mu )$?
One sufficient condition can be formulated with the help of the following
notion.

\begin{definition}
Let $f:\Omega \rightarrow X$ be an RL-integrable function. We say that a
function $g:\Omega \rightarrow X$ is $f$'s \emph{Bochner-integrable
equivalent} if $g$ is Bochner-integrable and for any measurable $A\subset
\Omega $,

\[
(RL)\int_Afd\mu =(Bochner)\int_Agd\mu \text{.} 
\]
\end{definition}

It is easy to see that if every $X$-valued RL-integrable function has a
Bochner-integrable equivalent, then $L_1(X,\mu )$ is complemented in $%
RL_1(X,\mu )$. Indeed, the equivalent Bochner-integrable function for the
given RL-integrable function $f$ is $f$'s projection onto $L_1(X,\mu )$. In
view of this, let us investigate under what conditions on the space $X$
every RL-integrable function has a Bochner-integrable equivalent.

One easy sufficient condition for $X$ is to have the Radon--Nikod\'ym
property (RNP).

\begin{theorem}
Let the Banach space $X$ have the RNP. Then any RL-integrable function $%
f:\Omega \rightarrow X$ has a Bochner-integrable equivalent.
\end{theorem}

\emph{Proof. }Consider the vector measure $\nu :\Sigma \rightarrow X$
defined by $\nu (A)=\int_Afd\mu $. It is countably-additive and has bounded
variation (this follows from the existence of integrable majorant of $f$).
Since $X$ possesses the RNP there must exist a Bochner-integrable function $%
g:\Omega \rightarrow X$ such that $\nu (A)=\int_Agd\mu $ for any measurable $%
A$. Obviously, $g$ is $f$'s Bochner-integrable equivalent. $\Box $

\bigskip Another simple condition is the following:

\begin{theorem}
\label{th-sep-compl-Bochner-equiv}Let $f:\Omega \rightarrow X$ be an
RL-integrable function and the set $\{\int_Afd\mu :\,A\in \Sigma \}$ be
contained in a separable complemented subspace $Y$ of $X$. Then $f$ has a
Bochner-integrable equivalent.
\end{theorem}

\emph{Proof. }Let $P$ be the projection from $X$ onto $Y$. Put $g=Pf$. Since 
$g$ is an image of $f$ under continuous linear map, it is RL-integrable.
Moreover, the values of $g$ lie in the separable space $Y$. Hence $g$ is
also Bochner integrable. Obviously, $\int_Agd\mu =P\int_Afd\mu $ for any $%
A\in \Sigma $, and since all integrals $\int_Afd\mu $ lie in $Y$, $%
P\int_Afd\mu =\int_Afd\mu $. Thus $g$ is the Bochner-integrable equivalent
of $f$. $\Box $

\bigskip Now let us note the following fact.

\begin{theorem}
\label{th-RL-integrals-sep}Let $f:\Omega \rightarrow X$ be an RL-integrable
function. Then the set $\{\int_Afd\mu :\,A\in \Sigma \}$ is separable.
\end{theorem}

\emph{Proof. }Fix a sequence of positive numbers $\epsilon _n\rightarrow 0$.
Since $f$ is RL-integrable, for each $n$ there exists a partition $\Pi _n$
such that for any two finer partitions $\Gamma ^{\prime }$ and $\Gamma
^{\prime \prime }$ with any sets of sampling points $T^{\prime }$ and $%
T^{\prime \prime }$ we have 
\begin{equation}
\left\| S(f,\Gamma ^{\prime },T^{\prime })-S(f,\Gamma ^{\prime \prime
},T^{\prime \prime })\right\| <\epsilon _n.  \label{eqn-Bochner-equiv-sep-1}
\end{equation}
For each $\Pi _n$ fix a set of sampling points $T_n$. Obviously, $S(f,\Pi
_n,T_n)\rightarrow \int fd\mu $. Let $U=\bigcup_{n=1}^\infty T_n$. Since $U$
is countable, the subspace $Y=\overline{\limfunc{Lin}}(f(U))$ is separable.
Obviously, $S(f,\Pi _n,T_n)\in \limfunc{Lin}(f(U))$ and therefore $\int
fd\mu \in Y$. Let us now show that in fact $\int_Afd\mu \in Y$ for all $A\in
\Sigma $.

Indeed, take an arbitrary $A\in \Sigma $. Since $U$ is negligible, $%
\int_Afd\mu =\int_{A\cup U}fd\mu $. Therefore, without loss of generality we
may and do assume that $U\subset A$. Denote by $\Pi _n^A$ the partition of $%
A $ formed by intersection of subsets of $\Pi _n$ with $A$. Then for any two
partitions of $A$, $\Gamma ^{\prime }\succ \Pi _n^A$ and $\Gamma ^{\prime
\prime }\succ \Pi _n^A$ with any sets of sampling points $T^{\prime }$ and $%
T^{\prime \prime }$ respectively, condition~(\ref{eqn-Bochner-equiv-sep-1})
remains true. Therefore $S(f,\Pi _n^A,T_n)\rightarrow \int_Afd\mu $ and
hence, $\int_Afd\mu \in Y.$ $\Box $

\bigskip Theorems~\ref{th-sep-compl-Bochner-equiv} and~\ref
{th-RL-integrals-sep} together give us the following useful corollary:

\begin{theorem}
Let the Banach space $X$ have the following property: every separable
subspace of $X$ is contained in a separable complemented subspace of $X$
(the class of such spaces includes, for example, all the WCG spaces). Then
any RL-integrable function $f:\Omega \rightarrow X$ has a Bochner-integrable
equivalent.
\end{theorem}

Thus, we have shown that in quite a wide class of Banach spaces every
RL-integrable function has a Bochner-integrable equivalent. However, there
exist spaces where this is not true. Let us show an example.

\emph{Example. }Consider the function $f:[0;1]\rightarrow L_\infty [0;1]$
defined by $f(t)=\chi _{[t;1]}$. We prove that it is RL-integrable, but does
not have an equivalent Bochner-integrable function.

Indeed, let us show the following identity: 
\begin{equation}
(\int_Afd\mu )(t)=\mu (A\cap [0;t])  \label{eqn-nonequiv-star}
\end{equation}
for any $t\in [0;1]$ and any Borel set $A$.

Since both sides of the equality contain $L_\infty $-valued
countably-additive measures, it suffices to show the statement for $A=[\frac
k{2^n};\frac{k+1}{2^n}]$, where $n\in \Bbb{N}\cup \{0\}$, $0\leq k\leq 2^n$.

To this end, we take an arbitrary $\epsilon >0$ and find such a positive
integer $N$ that $\frac 1{2^nN}<\epsilon $. Let us partition $A$ into $N$
equal intervals $\{\Delta _i\}_{i=1}^N$ of length $\frac 1{2^nN}$ and let $%
\{\Delta _{ij}\}_{i=1}^N\,_{j=1}^{n_i}$ be an arbitrary finer partition into
non-empty intervals, where $\bigcup_{j=1}^{n_i}\Delta _{ij}=\Delta _i$, $%
i=1,\ldots ,N$. Let also $\{t_{ij}\}$ be a set of sampling points for this
finer partition. We show that for the functions 
\[
S(t)=\left\langle \sum_{i=1}^N\sum_{j=1}^{n_i}f(t_{ij})\mu (\Delta
_{ij}),t\right\rangle 
\]
and 
\[
g(t)=\mu (A\cap [0;t]) 
\]
the condition 
\[
\left| S(t)-g(t)\right| <\epsilon 
\]
holds for any $t\in [0;1]$.

Indeed, for $t\in [0;\frac k{2^n}]$ we have 
\begin{equation}
S(t)=0=g(t).  \label{eqn-nonequiv-1}
\end{equation}
If $t\in [\frac{k+1}{2^n};1],$ then 
\begin{equation}
S(t)=\mu (A)=g(t)  \label{eqn-nonequiv-2}
\end{equation}
For $t\in \Delta _{i_0j_0}$ we can estimate: 
\[
S(t)=\sum_{i=1}^{i_0}\sum_{j=1}^{j_0}\mu (\Delta _{ij})\leq
\sum_{i=1}^{i_0}\sum_{j=1}^{n_{i_0}}\mu (\Delta _{ij})=\sum_{i=1}^{i_0}\mu
(\Delta _i)=\frac{i_0}{2^nN}, 
\]

\[
S(t)\geq \sum_{i=1}^{i_0}\sum_{j=1}^{j_0-1}\mu (\Delta _{ij})\geq
\sum_{i=1}^{i_0-1}\sum_{j=1}^{n_{i_0-1}}\mu (\Delta
_{ij})=\sum_{i=1}^{i_0-1}\mu (\Delta _i)=\frac{i_0-1}{2^nN}\text{, for }%
i_0>1, 
\]

\[
S(t)\geq 0\text{, for }i_0=1. 
\]

Analogously, the same estimates are checked for the function $g$.Thus, for $%
t\in \Delta _{i_0j_0}$ we have 
\begin{equation}
\left| S(t)-g(t)\right| \leq \frac 1{2^nN}<\epsilon .  \label{eqn-nonequiv-3}
\end{equation}
Combining~(\ref{eqn-nonequiv-1}), (\ref{eqn-nonequiv-2}) and (\ref
{eqn-nonequiv-3}), we obtain 
\[
\Vert S-g\Vert _\infty <\epsilon . 
\]

Since $\epsilon $, $\{\Delta _{ij}\}$ and $\{t_{ij}\}$ have been chosen
arbitrarily, equality~(\ref{eqn-nonequiv-star}) is proved, meaning by $%
\int_Afd\mu $ the Riemann integral (for $A$ of the form $[\frac k{2^n};\frac{%
k+1}{2^n}]$). But since the Riemann integrability implies RL-integrability,
we have shown a stronger result.

It remains to prove that $f$ does not allow a Bochner-integrable equivalent.
Since the Bochner integral is differentiable as the function of the upper
limit, it suffices to show that the function 
\[
G(t)=\int_0^tfd\mu 
\]
is not differentiable at any point $t\in (0;1)$.

Let $t_0\in (0;1)$, $\Delta t>0$. Then 
\[
\frac{G(t_0+\Delta t)-G(t_0)}{\Delta t}=\frac 1{\Delta
t}\int_{t_0}^{t_0+\Delta t}fd\mu . 
\]

But equality~(\ref{eqn-nonequiv-star}) shows that the functions $\frac
1{\Delta t}\int_{t_0}^{t_0+\Delta t}fd\mu $ do not form a fundamental
family, as $\Delta t\rightarrow 0$. So, this family has no limit and hence $%
G(t)$ is not differentiable.

\bigskip In the example above, note that the values of $f$ are functions
with at most one discontinuity of the first order. The space of all such
functions is isomorphic to a subspace of the space $C(K)$, where $K$ is the
topological space known as ``two arrows of Alexandrov''. On the other hand,
the values of $\int_Afd\mu $ are all contained in $C[0;1]$. Since $f$ has
been shown to have no a Bochner-integrable equivalent, theorem~\ref
{th-sep-compl-Bochner-equiv} implies the following interesting corollary:

\begin{corollary}
The space $C[0;1]$ is not complemented in $C$ on ``two arrows'' and is not
contained in any separable complemented subspace thereof.
\end{corollary}

Another question that seems natural to ask is under what condition $%
RL_1(X,\mu )$ coincides with $L_1(X,\mu )$, i.e.~when every RL-integrable
function is also Bochner-integrable. Let us present one sufficient condition.

\begin{theorem}
Let $X$ be a Banach space such that every $X$-valued RL-integrable function
has a Bochner-integrable equivalent. Let $X$ have a countable set of
functionals separating the points of $X$. Then every RL-integrable $X$%
-valued function is Bochner-integrable.
\end{theorem}

\emph{Proof.} Let $f:\Omega \rightarrow X$ be an RL-integrable function and $%
g\in L_1(\Omega ,X)$ its Bochner-integrable equivalent. We show that $f=g$
almost everywhere. Indeed, for any $x^{*}\in X^{*}$ and any measurable $%
A\subset \Omega $ consider: 
\begin{eqnarray*}
\int_Ax^{*}(f-g)d\mu &=&\int_Ax^{*}fd\mu -\int_Ax^{*}gd\mu = \\
&=&x^{*}(\int_Afd\mu )-x^{*}(\int_Agd\mu )=0,
\end{eqnarray*}
where $\int_Afd\mu $ denotes $f$'s RL-integral and $\int_Agd\mu $ denotes $g$%
's Bochner integral. These two integrals are equal, since $g$ is $f$'s
Bochner-integrable equivalent. Thus the Lebesgue integral of $x^{*}(f-g)$ is
zero over any measurable set $A$, which means that $x^{*}(f-g)=0$ a.e.

Now let $\{x_n^{*}\}_{n=1}^\infty $ be the countable set of functionals
separating the points of $X$. Put $A_n=\{t\in \Omega :\,x_n^{*}(f-g)\neq 0\}$
and let $A=\bigcup_{n=1}^\infty A_n$. Note that for any $t\in \Omega
\setminus A$ and any functional $x_n^{*}$, $x_n^{*}(f-g)(t)=0$. Since $A$ is
negligible, $x_n^{*}(f-g)=0$ a.e., therefore $f$ itself is
Bochner-integrable. $\Box $

\bigskip\ Let us further investigate the isomorphic structure of $RL_1(X,\mu
)$. We show that this space can be ``very large'', meaning it can contain an
isomorphic copy of the space $l_\infty (\Gamma )$.

\emph{Example.} Suppose $X=l_2([0,1])$ and $\mu $ is the Lebesgue measure.
Then $RL_1(X,\mu )$ contains an isometrical copy of $l_\infty ([0,1])$.

First note that there is a continuum-cardinality family of non-measurable
mutually disjoint sets $\{A_t\}_{t\in [0,1]}$, where $A_t\subset [0;1]$ and
the outer measure of each $A_t$ equals $1$. Such a construction can be
found, for example, in \cite{Kadets-Tseytlin}, in the proof of theorem~2.16.

Now consider the following function $f:[0,1]\rightarrow l_2([0,1])$: $%
f(t)=e_t$, where $e_t$ is an ort of $l_2([0,1])$: $e_t=\chi _{\{t\}}$.
Define the linear map $U:l_\infty ([0,1])\rightarrow RL_1(X,\mu )$ as
follows: 
\[
U(\alpha )=\sum_{t\in [0,1]}\alpha _tf\cdot \chi _{A_t}\text{,} 
\]
whenever $\alpha =(\alpha _t)_{t\in [0,1]}\in l_\infty ([0,1])$.

To see that $U(\alpha )$ is indeed in $RL_1(X,\mu )$ for every $\alpha \in
l_\infty ([0,1])$, take any partition $\Pi =\{\Delta _i\}$ and any sample
points $U=\{t_i\}$ and estimate: 
\begin{eqnarray*}
\left\| S(U(\alpha ),\Pi ,T)\right\| &\leq &\Vert \alpha \Vert _\infty
\left\| \stackrel{\infty }{\sum_{i=1}}e_{t_i}\mu (\Delta _i)\right\| = \\
&=&(\stackrel{\infty }{\sum_{i=1}}\mu (\Delta _i)^2)^{1/2}\leq (\stackrel{%
\infty }{\sum_{i=1}}d(\Pi )\cdot \mu (\Delta _i))^{1/2}= \\
&=&\sqrt{d(\Pi )}(\stackrel{\infty }{\sum_{i=1}}\mu (\Delta _i))^{1/2}\leq 
\sqrt{d(\Pi )}(\sum_{i=1}^\infty \mu (\Delta _i))^{1/2}=\sqrt{d(\Pi )},
\end{eqnarray*}
where $d(\Pi )$ denotes $\sup \{\mu (\Delta _i)\}$. Now take an arbitrary $%
\epsilon >0$. Fix a partition $\Pi $ with $d(\Pi )<\epsilon $. Then for any $%
\Pi ^{\prime }\succ \Pi $, any $T^{\prime }$, $\left\| S(U(\alpha ),\Pi
^{\prime },T^{\prime })\right\| \leq \sqrt{d(\Pi ^{\prime })}\leq \sqrt{%
d(\Pi )}\leq \epsilon $. Thus $U(\alpha )$ is RL integrable and the integral
equals $0$.

Further, it is clear that $\Vert U(\alpha )(t)\Vert \leq \Vert \alpha \Vert
_\infty $ for every $t\in [0,1]$. So $\Vert U(\alpha )\Vert _{RL}\leq \Vert
\alpha \Vert _\infty $. On the other hand, for a fixed $n\in \Bbb{N}$ any
integrable majorant of $\Vert U(\alpha )(t)\Vert $ must be not less than $%
\Vert \alpha \Vert _\infty -\frac 1n$ on a set of full outer measure, and
therefore almost everywhere. Therefore $\Vert U(\alpha )\Vert _{RL}\geq
\Vert \alpha \Vert _\infty -\frac 1n$, $n\in \Bbb{N}$ $.$ So $\Vert U(\alpha
)\Vert _{RL}=\Vert \alpha \Vert _\infty $ and we are done.

\section{Limit Set $I(f)$}

Even if a function is not RL-integrable, we still can consider the set of
limit points of its RL-integral sums, which, in a sense, plays the role of
RL-integral.

\begin{definition}
Let $f:\Omega \rightarrow X$ be an arbitrary function with values in a
Banach space $X$. We say that a point $x\in X$ belongs to the limit set $%
I(f) $ if for every $\epsilon >0$ and any partition $\Pi $ there exists a
partition $\Gamma \succ \Pi $ and a set of sampling points $T$ such that $%
||S(f,\Gamma ,T)-x||<\epsilon $, $S(f,\Gamma ,T)$ being an absolute integral
sum. In other words, $I(f)$ is the set of all limit points of the net of $f$%
's absolute integral sums.
\end{definition}

Let us point out two important facts about the set $I(f)$ (the proofs can be
found in \cite{Kadets-Tseytlin}):

\begin{enumerate}
\item  The limit set $I(f)$ is always convex.

\item  For an arbitrary convex closed set $S$ of no more than continuum
cardinality in a Banach space $X$, there exists a function $%
 f:[0;1]\rightarrow X$ such that $I(f)=S$.
\end{enumerate}

We will show below that the problem of existence of such limit sets is not
always solved in positive. However, as the following theorem shows, for a
large class of Banach spaces it is.

\begin{theorem}
\label{th-WCG-nonempty}Let $(\Omega ,\Sigma ,\mu )$ be a measure space, $X$
a WCG space and $f\in \overline{L_1}(\Omega ,\Sigma ,\mu ,X)$ a function.
Then $I(f)$ is not empty.
\end{theorem}

The proof is almost entirely contained in the following lemmas.

\begin{lemma}
\label{lem-WCG-1}Let $(\Omega ,\Sigma ,\mu )$ be a probability measure
space, $X$ a WCG space generated by a convex balanced weakly-compact set $K$%
. Suppose $f:\Omega \rightarrow X$ is an arbitrary function, not necessarily
measurable. Then for any sequence of positive real numbers $\epsilon _k$
there exists a sequence of functions $g_k:\Omega \rightarrow X$, sequence of
subsets $A_k\subset \Omega $ and sequence of integers $n_k\in \Bbb{N}$ so
that the following properties hold:

\begin{enumerate}
\item  $\mu ^{*}(\bigcap_{i=k}^NA_i)>1-\sum_{i=k}^N\epsilon _i$ for any $%
k\leq N$;

\item  $g_k(A_k)\subset n_kK$;

\item  $\left\| g_k(t)-f(t)\right\| \leq \epsilon _k$ for $t\in A_k$ and $%
g_k(t)=0$ for $t\in \Omega \setminus A_k$.
\end{enumerate}
\end{lemma}

\emph{Proof.} Note that $\bigcup_{n=1}^\infty (nK+\epsilon _kB(X))=X$ for
any $k$. Due to this, for any $\epsilon _k$ we can choose an index $n$ so
large that $\mu ^{*}(f^{-1}(nK+\epsilon _kB(X)))>1-\epsilon _k$. Moreover,
for an arbitrary $B\subset \Omega $ we can choose $n$ so large that $\mu
^{*}((f^{-1}(nK+\epsilon _kB(X)))\cap B)\geq \mu ^{*}(B)-\epsilon _k$. These
two observations allow us to construct the sequence $\{n_k\}\subset \Bbb{N}$
by induction, so that if $A_k=f^{-1}(n_kK+\epsilon _kB(X))$, then

\[
\mu ^{*}(A_k)>1-\epsilon _k\text{,} 
\]
for any $k$, and 
\[
\mu ^{*}(A_k\cap (\bigcap_{i=j}^{k-1}A_i))\geq \mu
^{*}(\bigcap_{i=j}^{k-1}A_i)-\epsilon _k\text{,} 
\]
for any $k$ and any $j<k$.

It is easy to see now that the first of lemma's conditions is satisfied: 
\begin{eqnarray*}
\mu ^{*}(\bigcap_{i=k}^NA_i) &=&\mu ^{*}(A_N\cap
(\bigcap_{i=k}^{N-1}A_i))\geq \mu ^{*}(\bigcap_{i=k}^{N-1}A_i)-\epsilon _N \\
&\geq &\mu ^{*}(\bigcap_{i=k}^{N-2}A_i)-\epsilon _N-\epsilon _{N-1}\geq
\ldots \geq \mu (A_k)-\sum_{i=k+1}^N\epsilon _i \\
&\geq &1-\sum_{i=k}^N\epsilon _i\text{.}
\end{eqnarray*}
Now note that by the construction, 
\begin{equation}
f(A_k)\subset n_kK+\epsilon _kB(X)\text{.}  \label{eq-WCG-1}
\end{equation}
This allows us to construct functions $g_k$ satisfying conditions 2 and 3 of
the lemma. Indeed, we put $g_k(t)=0$ for $t\in \Omega \setminus A_k$. If $%
t\in A_k,$ then by~(\ref{eq-WCG-1}), there exists a point $y\in n_kK$ such
that $\Vert f(t)-y\Vert \leq \epsilon _k$. We put then $g_k(t)=y$. This
finishes the proof. $\Box $

\begin{definition}
Let $g:\Omega \rightarrow X$ be a function, $A\subset \Omega $. We denote by 
$\sigma (g|A)$ the set of all integral sums of $g$, \emph{controlled by} $A$%
, that is, a sampling point is always chosen in $A$, whenever the partition
subset intersects with $A.$ We also denote by $I(g|A)$ the set of limits of
the integral sums from $\sigma (g|A)$.
\end{definition}

\begin{lemma}
\label{lem-WCG-2}Suppose a function $f:\Omega \rightarrow X$ has an
integrable majorant $h$: $\Vert f(t)\Vert \leq h(t)$, $h\in L_1(\Omega
,\Sigma ,\mu )$. Let the real numbers $\epsilon _n>0$ be so small that if $%
A\in \Sigma $, $\mu (A)<\epsilon _n$, then 
\begin{equation}
\int_Ahd\mu <2^{-n-1}\text{.}  \label{eqn-WCG-2}
\end{equation}
Assume also that $\epsilon _1<\frac 14$ and $\epsilon _{n+1}<\frac
1{2^n}\epsilon _n$. Then, under the conditions of lemma~\ref{lem-WCG-1},
there exists a convergent sequence $\{x_n\}\subset X$ such that $x_n\in
I(g_n|A_n)$.
\end{lemma}

\emph{Proof.} On the set of all partitions define an ultrafilter $\mathcal{U}
$, which dominates the filter of the refinement direction.

For every partition $\Gamma =\{\Delta _i\}_{i=1}^\infty $ define a number $%
N(\Gamma )\in \Bbb{N}$ so that $N(\Gamma )\rightarrow \infty $, as $\Gamma $
gets finer (e.g.,~$N(\Gamma )=[1/\sup_i\mu (\Delta _i)]$). Construct the
integral sums 
\[
S_k(\Gamma )=\sum_{j=1}^\infty g_k(t_j^{(k)}(\Gamma ))\mu (\Delta _j)\in
\sigma (g_k|A_k)\text{,} 
\]
where $k=1,2,\ldots ,N(\Gamma )$, so that the following condition is
satisfied: if for some $j\in \Bbb{N}$ there exists an index $s\leq N(\Gamma
) $ such that 
\[
\Delta _j\cap \bigcap_{i=s}^{N(\Gamma )}A_i\neq \emptyset 
\]
(denote the smallest such index by $s(j)$), then 
\[
t_j^{(k)}(\Gamma )=t_j^{(s(j))}(\Gamma )\in \bigcap_{i=s(j)}^{N(\Gamma )}A_i 
\]
for $k=s(j),s(j)+1,\ldots ,N(\Gamma )$.

Note that under such a choice of sampling points, in view of the first
condition of lemma~\ref{lem-WCG-1}, for each $k$ the total measure of all
those $\Delta _j$ where $t_j^{(k)}=t_j^{(k+1)}$, is bounded below by 
\[
1-\sum_{i=k}^{N(\Gamma )}\epsilon _k>1-\epsilon _{k-1}\text{.} 
\]

Since $h+\epsilon _k$ is a majorant for $g_k$, we have 
\begin{equation}
\lim_{\mathcal{U}}\left\| S_k(\Gamma )-S_{k+1}(\Gamma )\right\| \leq
2\epsilon _k+2^{-k}\text{.}  \label{eqn-WCG-3}
\end{equation}
Since for any $k$, given a sufficiently fine partition $\Gamma $, the number 
$N(\Gamma )$ is greater than $k$, the sum $S_k(\Gamma )$ will be correctly
defined and, due to condition 2 of lemma~\ref{lem-WCG-1}, $S_k(\Gamma )\in
n_kK$. From the compactness argument we infer that there exists a weak limit 
$w-\lim_{\mathcal{U}}S_k(\Gamma )$. Denote this limit by $x_k$. It was
proved in \cite{Kadets-Tseytlin} that actually $x_k\in I(g_k|A_k)$. Due to (%
\ref{eqn-WCG-3}), $\{x_k\}$ forms a fundamental sequence. Thus $\{x_k\}$ is
required. $\Box $

\bigskip \emph{Proof of theorem~\ref{th-WCG-nonempty}.}

Let $\{x_k\}$ be the sequence from lemma~\ref{lem-WCG-2}, $%
x=\lim_{k\rightarrow \infty }x_k$. Due to the condition $\mu
^{*}(A_k)>1-\epsilon _k$ , condition 3 of lemma~\ref{lem-WCG-1} and (\ref
{eqn-WCG-2}), each $x_k$ can be approximated with the precision of $\epsilon
_k+\frac 1{2^k}$ by arbitrary fine integral sums of $f$. Therefore $x\in
I(f) $, which is to be proved. $\Box $

\bigskip So, we have shown that $I(f)$ is non-empty for any $f\in \overline{%
L_1}(X)$ for quite a large class of WCG spaces. This result supersedes the
one of \cite{Kadets-Tseytlin}, where an analogous theorem is shown for
separable and reflexive spaces $X$. However, there exist spaces where $I(f)$
can be empty for certain functions $f\in \overline{L_1}(X)$. Let us exhibit
an example.

\emph{Example. }Consider the space $X=l_1([0;1])$ and function $%
f:[0;1]\rightarrow l_1([0;1])$ given by $f(t)=e_t$, i.e.~$f(t)=\chi _{\{t\}}$%
. Obviously, $f$ is bounded and therefore $f\in \overline{L_1}(X)$. It is
easy to see that the $l_1$-norm of any integral sum $S(f,\Gamma ,T)$ of $f$
equals $1$. Suppose that there exists an $x\in I(f)$. Then $||x||_{l_1}=1$.
Consider the coordinate functionals $\delta _t\in (l_1[0;1])^{*}$: $\delta
_t(g)=g(t)$. For every $t\in [0;1]$ we have $\delta _t(f(\tau ))=0$ for
almost all $\tau \in [0;1]$ (in fact, for all $\tau \neq t$). Since $x\neq 0$%
, there exists a $t_0\in [0;1]$ such that $\delta _{t_0}(x)\neq 0$. On the
other hand, 
\[
\delta _{t_0}(x)\in \delta _{t_0}(I(f))\subset I(\delta _{t_0}\circ f)=\{0\}%
\text{.} 
\]
This contradiction shows that $I(f)=\oslash $.

\bigskip If a function $f$ is RL-integrable, then $I(f)$ obviously consists
of a single point, $f$'s integral. The converse is not always true: a
function may have a single-point $I(f)$ and still not be RL-integrable. The
following weaker statement is true, however.

\begin{theorem}
\label{th-single-point-Pettis}Let $X$ be a Banach space such that every
function $f\in \overline{L_1}(\Omega ,\Sigma ,\mu ,X)$ from any measure
space $(\Omega ,\Sigma ,\mu )$ has a non-empty limit set $I(f)$. Let a
function $f\in \overline{L_1}(\Omega ,\Sigma ,\mu ,X)$ be such that $I(f)$
consists of a single point $x$. Then $f$ is Pettis integrable and $x$ is $f$%
's Pettis integral.
\end{theorem}

Let us prove some lemmas first.

\begin{lemma}
\label{lem-single-point-restriction}Let $X$ and $f$ be such as in theorem~%
\ref{th-single-point-Pettis} and let $A\in \Sigma $. Then for the
restriction $f|_A$ of $f$ on $A$, $I(f|_A)$ is a singleton.
\end{lemma}

\emph{Proof.} Due to the properties of $X$, $I(f|_A)$ cannot be empty. To
prove that $I(f|_A)$ is a singleton, assume there are two distinct points $%
x_1,x_2\in I(f|_A)$. Fix an $\epsilon >0$. Denote $B=\Omega \setminus A$ and
pick any $y\in I(f|_B)$.

Now take an arbitrary partition $\Pi $ of $\Omega $. Let $\Pi ^{\prime }$ be
a partition that is finer than both $\Pi $ and the partition of $\Omega $
into two sets $A$ and $B$. Then every member set of $\Pi ^{\prime }$ is
either a subset of $A$ or a subset of $B$. Therefore we can consider the
partition $\Pi ^A$ of $A$ formed by those members of $\Pi ^{\prime }$ that
lie within $A$, and a partition $\Pi ^B$ of $B$ formed by those members of $%
\Pi ^{\prime }$ that lie within $B$.

Since $x_1\in I(f|_A)$, there exists a partition $\Pi _1^A\succ \Pi ^A$ and
a set of sampling points $T_1^A$ for it such that $\left\| S(f|_A,\Pi
_1^A,T_1^A)-x_1\right\| <\epsilon /2$. Since also $x_2\in I(f|_A)$, there
exists a partition $\Pi _2^A\succ \Pi ^A$ and a set of sampling points $%
T_2^A $ for it such that $\left\| S(f|_A,\Pi _2^A,T_2^A)-x_2\right\|
<\epsilon /2$. And since $y\in I(f|_B)$, there exists a partition $\Pi
_1^B\succ \Pi ^B$ and a set of sampling points $T_1^B$ for it such that $%
\left\| S(f|_B,\Pi _1^B,T_1^B)-y\right\| <\epsilon /2$. Combine the
partitions $\Pi _1^A$ and $\Pi _1^B$ into partition $\Pi _1$ of the entire $%
\Omega $, and put $T_1=T_1^A\cup T_1^B$. Then $\left\| S(f,\Pi
_1,T_1)-(x_1+y)\right\| <\epsilon $. At the same time combine the partitions 
$\Pi _2^A$ and $\Pi _1^B $ into partition $\Pi _2$ of the entire $\Omega $,
and put $T_2=T_2^A\cup T_1^B$. Then $\left\| S(f,\Pi _2,T_2)-(x_2+y)\right\|
<\epsilon $. Since $\Pi _1\succ \Pi $ and $\Pi _2\succ \Pi $, both $x_1+y$
and $x_2+y$ belong to $I(f)$, which is impossible. Hence $I(f|_A)$ consists
of a single point. $\Box $

\begin{lemma}
Let $X$ and $f$ be such as in theorem~\ref{th-single-point-Pettis}. Then $f$
is weakly measurable.
\end{lemma}

\emph{Proof.} Assume the contrary. Then there exists a functional $x^{*}\in
X^{*}$ such that $x^{*}f$ is a non-measurable function. Since $f\in 
\overline{L_1}(\Omega ,\Sigma ,\mu ,X)$, $x^{*}f$ must have an integrable
(and hence measurable) majorant. Therefore, there must exist the smallest
measurable majorant $f_2$ of $x^{*}f$ and the largest measurable minorant $%
f_1$ of $x^{*}f$. In other words, if $g:\Omega \rightarrow \Bbb{R}$ is
measurable and $x^{*}f\leq g$ a.e., then $f_2\leq g$ a.e., and similarly, if 
$g:\Omega \rightarrow \Bbb{R}$ is measurable and $x^{*}f\geq g$ a.e., then $%
f_1\geq g$ a.e. Note that $f_1$ and $f_2$ cannot coincide almost everywhere,
since that would mean that $f_1=x^{*}f=f_2$ a.e. and $x^{*}f$ would be
measurable. Note that $\{t:\,f_1(t)\neq f_2(t)\}=\bigcup_{n=1}^\infty
\{t:\,f_2(t)-f_1(t)>1/n\}$. Since the set at the left-hand side is
non-negligible, one of the sets on the right-hand side must be
non-negligible too. Therefore, there exists a non-negligible measurable set $%
A$ and an $\epsilon >0$ such that $f_1(t)<f_2(t)-\epsilon $ for any $t\in A$.

Consider $f|_A$. Due to lemma~\ref{lem-single-point-restriction}, $I(f|_A)$
consists of a single point. Consider the following two sets: 
\begin{eqnarray*}
A_1 &=&\{t\in A:\,x^{*}f(t)>\frac 23f_2(t)+\frac 13f_1(t)\} \\
A_2 &=&\{t\in A:\,x^{*}f(t)<\frac 23f_1(t)+\frac 13f_2(t)\}
\end{eqnarray*}
and let $B_1=A\setminus A_1$, $B_2=A\setminus A_2$. Note that $B_1$ cannot
contain any measurable non-negligible set. Indeed, assume that $C$ is a
measurable set, $\mu (C)>0$ and $C\subset B_1$. This means that for any $%
t\in C$, $x^{*}f(t)\leq \frac 23f_2(t)+\frac 13f_1(t)<f_2(t)$. But now
consider function $g$, which is equal to $\frac 23f_2(t)+\frac 13f_1(t)$ for 
$t\in C$ and coincides with $f_2$ outside of $C$. This function is
measurable, it is a majorant of $x^{*}f$, but $g<f$ on a non-negligible set $%
C$. This contradicts the definition of $f_2$ as the smallest measurable
majorant of $x^{*}f$. Thus we have shown that $B_1$ contains no measurable
non-negligible subset, which means that $\mu _{*}(B_1)=0$ and $\mu
^{*}(A_1)=\mu (A)$. It is easy to apply the same argument to $B_2$ and $A_2$
to show that $\mu ^{*}(A_2)=\mu (A)$.

Now consider a $\sigma $-field $\Sigma _{A_1}$ of subsets of $A_1$ of the
form $C\cap A_1$, where $C\in \Sigma $. Define $\mu |_A(C\cap A_1)=\mu
(C\cap A)$. So, we obtain a measure space $(A_1,\Sigma _{A_1},\mu |_A)$. It
is easy to verify that this space is correctly defined, since $\mu
^{*}(A_1)=\mu (A)$. The restriction $f|_{A_1}$ is a function from this
measure space to the Banach space $X$. Due to the properties of $X$, there
exists an $x_1\in I(f|_{A_1})$. By an analogous argument we can construct
the measure space $(A_2,\Sigma _{A_2},\mu |_A)$ and find a point $x_2\in
I(f|_{A_2})$.

Let us show that $x_1\neq x_2$. Indeed, note that $x^{*}(x_1)\in
I(x^{*}f|_{A_1})$ and $x^{*}(x_2)\in I(x^{*}f|_{A_2})$. Consider an integral
sum of $x^{*}f|_{A_1}$. It has the form $\sum x^{*}f(t_i)\mu (\Delta _i\cap
A_1)=\sum x^{*}f(t_i)\mu (\Delta _i)$, where $\{\Delta _i\}$ is a partition
of $A$ and $t_i\in \Delta _i\cap A_1$. Thus all integral sums of $%
x^{*}f|_{A_1}$ dominate the integral sums of the function $\frac
23f_2(t)+\frac 13f_1(t)$ over $A$. On the other hand, the same argument
shows that all integral sums of $x^{*}f|_{A_2}$ are dominated by the
integral sums of the function $\frac 23f_1(t)+\frac 13f_2(t)$ over $A$.
Since the values if these two functions differ by at least $\epsilon /3$ at
all points of $A$, this implies that $x^{*}(x_1)>x^{*}(x_2)$, which means
that $x_1\neq x_2$.

Let us show that $x_1,x_2\in I(f|_A)$. Indeed, since $\mu ^{*}(A_1)=\mu
^{*}(A_2)=\mu (A)$, any integral sum over $A_1$ of the form $\sum f(t_i)\mu
(\Delta _i\cap A_1)$ is equal to the integral sum $\sum f(t_i)\mu (\Delta
_i) $ over $A$, and the same is true for $A_2$. So, we have found two
different points $x_1$ and $x_2$ in $I(f|_A)$, which is impossible. This
contradiction proves that $f$ is weakly measurable. $\Box $

\bigskip\emph{Proof of theorem~\ref{th-single-point-Pettis}. }

The previous lemma shows that the function $f$ is weakly measurable. Take
any measurable $A\subset \Omega $. Take an $x^{*}\in X^{*}$. The real-valued
function $x^{*}f|_A$ is measurable. Since $f$ has an integrable majorant, so
does $x^{*}f|_A$. Therefore $x^{*}f|_A$ is Lebesgue-integrable and thus,
RL-integrable. Hence $I(x^{*}f|_A)$ consists of a single point, $%
\int_Ax^{*}fd\mu $. Let $x_A$ be the only point of $I(f|_A)$. Since $%
x^{*}(x_A)\in I(x^{*}f|_A)$, we have $x^{*}(x_A)=\int_Ax^{*}fd\mu $, and
this holds for any functional $x^{*}\in X^{*}$ and any subset $A\in \Sigma $%
. So $f$ is Pettis integrable and $\int fd\mu =x$. $\Box $

\end{document}